\documentclass{amsart}
\usepackage{amssymb, amsbsy, amsthm, amsmath, amstext, amsopn, verbatim}
\usepackage[all]{xy}
\usepackage{amsfonts}
\usepackage{amscd}
\newtheorem{thm}{Theorem} [section]
\newtheorem{lemma}[thm]{Lemma}

\newtheorem{corollary}[thm]{Corollary}
\newtheorem{prop}[thm]{Proposition}
\newtheorem{notation}[thm]{Notation}
\newtheorem{convention}[thm]{Convention}

\theoremstyle{definition}

\newtheorem{defn}[thm]{Definition}

\newtheorem{example}[thm]{Example}

\theoremstyle{remark}

\newtheorem{remark}[thm]{Remark}

\begin{document}

\numberwithin{equation}{section}

\newcommand{\hs}{\mbox{\hspace{.4em}}}
\newcommand{\ds}{\displaystyle}
\newcommand{\bd}{\begin{displaymath}}
\newcommand{\ed}{\end{displaymath}}
\newcommand{\bcd}{\begin{CD}}
\newcommand{\ecd}{\end{CD}}

\newcommand{\proj}{\operatorname{Proj}}
\newcommand{\bproj}{\underline{\operatorname{Proj}}}
\newcommand{\spec}{\operatorname{Spec}}
\newcommand{\bspec}{\underline{\operatorname{Spec}}}
\newcommand{\bfspec}{\underline{\operatorname{FSpec}}}
\newcommand{\pline}{{\mathbf P} ^1}
\newcommand{\pplane}{{\mathbf P}^2}
\newcommand{\coker}{{\operatorname{coker}}}
\newcommand{\ldb}{[[}
\newcommand{\rdb}{]]}

\newcommand{\Sym}{\operatorname{Sym}^{\bullet}}
\newcommand{\Symp}{\operatorname{Sym}}
\newcommand{\Pic}{\operatorname{Pic}}
\newcommand{\AAut}{\operatorname{Aut}}
\newcommand{\PAut}{\operatorname{PAut}}

\newcommand{\good}{good\hspace{.15em}}
\newcommand{\too}{\twoheadrightarrow}
\newcommand{\C}{{\mathbf C}}
\newcommand{\cA}{{\mathcal A}}
\newcommand{\cS}{{\mathcal S}}
\newcommand{\cV}{{\mathcal V}}
\newcommand{\cM}{{\mathcal M}}
\newcommand{\bA}{{\mathbf A}}
\newcommand{\cB}{{\mathcal B}}
\newcommand{\cC}{{\mathcal C}}
\newcommand{\cD}{{\mathcal D}}
\newcommand{\D}{{\mathcal D}}
\newcommand{\cs}{{\mathbf C} ^*}
\newcommand{\boldc}{{\mathbf C}}
\newcommand{\cE}{{\mathcal E}}
\newcommand{\cF}{{\mathcal F}}
\newcommand{\cG}{{\mathcal G}}
\newcommand{\cH}{{\mathcal H}}
\newcommand{\h}{{\mathfrak h}}
\newcommand{\cJ}{{\mathcal J}}
\newcommand{\cK}{{\mathcal K}}
\newcommand{\cL}{{\mathcal L}}
\newcommand{\baL}{{\overline{\mathcal L}}}
\newcommand{\M}{{\mathcal M}}
\newcommand{\bM}{{\mathbf M}}
\newcommand{\cN}{{\mathcal N}}
\newcommand{\theo}{\mathcal{O}}
\newcommand{\cP}{{\mathcal P}}
\newcommand{\cR}{{\mathcal R}}
\newcommand{\bR}{{\mathbf R}}
\newcommand{\boldp}{{\mathbf P}}
\newcommand{\boldq}{{\mathbf Q}}
\newcommand{\bbL}{{\mathbf L}}
\newcommand{\cQ}{{\mathcal Q}}
\newcommand{\cO}{{\mathcal O}}
\newcommand{\Oo}{{\mathcal O}}
\newcommand{\OX}{{\Oo_X}}
\newcommand{\OY}{{\Oo_Y}}
\newcommand{\otY}{{\underset{\OY}{\ot}}}
\newcommand{\otX}{{\underset{\OX}{\ot}}}
\newcommand{\cU}{{\mathcal U}}\newcommand{\cX}{{\mathcal X}}
\newcommand{\cW}{{\mathcal W}}
\newcommand{\boldz}{{\mathbf Z}}
\newcommand{\qgr}{\operatorname{q-gr}}
\newcommand{\gr}{\operatorname{gr}}
\newcommand{\coh}{\operatorname{coh}}
\newcommand{\End}{\operatorname{End}}
\newcommand{\Hom}{\operatorname{Hom}}
\newcommand{\uHom}{\underline{\operatorname{Hom}}}
\newcommand{\uHomY}{\uHom_{\OY}}
\newcommand{\uHomX}{\uHom_{\OX}}
\newcommand{\Ext}{\operatorname{Ext}}
\newcommand{\uExt}{\underline{\operatorname{Ext}}}
\newcommand{\bExt}{\operatorname{\bf{Ext}}}
\newcommand{\Tor}{\operatorname{Tor}}

\newcommand{\inv}{^{-1}}
\newcommand{\airtilde}{\widetilde{\hspace{.5em}}}
\newcommand{\airhat}{\widehat{\hspace{.5em}}}
\newcommand{\nt}{^{\circ}}
\newcommand{\del}{\partial}

\newcommand{\supp}{\operatorname{supp}}
\newcommand{\GK}{\operatorname{GK-dim}}
\newcommand{\hd}{\operatorname{hd}}
\newcommand{\id}{\operatorname{id}}
\newcommand{\res}{\operatorname{res}}
\newcommand{\lrar}{\leadsto}
\newcommand{\im}{\operatorname{Im}}
\newcommand{\HH}{\operatorname{H}}
\newcommand{\TF}{\operatorname{TF}}
\newcommand{\Bun}{\operatorname{Bun}}
\newcommand{\Hilb}{\operatorname{Hilb}}
\newcommand{\nthord}{^{(n)}}
\newcommand{\Aut}{\underline{\operatorname{Aut}}}
\newcommand{\Gr}{\operatorname{\bf Gr}}
\newcommand{\Fr}{\operatorname{Fr}}
\newcommand{\GL}{\operatorname{GL}}
\newcommand{\SL}{\operatorname{SL}}
\newcommand{\ff}{\footnote}
\newcommand{\ot}{\otimes}
\def\Ext{\operatorname {Ext}}
\def\Hom{\operatorname {Hom}}
\def\Ind{\operatorname {Ind}}
\def\bbZ{{\mathbb Z}}

\newcommand{\nc}{\newcommand}
\newcommand{\on}{\operatorname}
\nc{\cont}{\on{cont}}
\nc{\rmod}{\on{mod}}
\nc{\Mtil}{\widetilde{M}}
\nc{\wb}{\overline}
\nc{\wt}{\widetilde}
\nc{\wh}{\widehat}
\nc{\mc}{\mathcal}
\nc{\Mbar}{\wb{M}}
\nc{\Nbar}{\wb{N}}
\nc{\Mhat}{\wh{M}}
\nc{\pihat}{\wh{\pi}}
\nc{\JYX}{\cJ_{Y\leftarrow X}}
\nc{\JXY}{\cJ_{X\rightarrow Y}}
\nc{\phitil}{\wt{\phi}}
\nc{\Qbar}{\wb{Q}}
\nc{\DYX}{\D_{Y\leftarrow X}}
\nc{\DXY}{\D_{X\to Y}}
\nc{\dR}{\stackrel{\bbL}{\underset{\D_X}{\ot}}}
\nc{\co}{\on{co}}
\nc{\aline}{{\mathbf A}^1}

\title{Cusps and $\cD$-Modules}
\author{David Ben-Zvi}
\address{Department of Mathematics\\University of Chicago\\Chicago, IL 60637}
\email{benzvi@math.uchicago.edu}
\author{Thomas Nevins}
\address{Department of Mathematics\\University of Michigan\\Ann Arbor, MI 
48109-1109}
\email{nevins@umich.edu}



\begin{abstract}
We study interactions between the categories of $\D$-modules on smooth
and singular varieties. For a large class of singular varieties $Y$,
we use an extension of the Grothendieck-Sato formula to show that
$\D_Y$-modules are equivalent to stratifications on $Y$, and as a
consequence are unaffected by a class of homeomorphisms, the {\em
cuspidal quotients}. In particular, when $Y$ has a smooth bijective
normalization $X$, we obtain a Morita equivalence of $\D_Y$ and $\D_X$
and a Kashiwara theorem for $\D_Y$, thereby solving conjectures of
Hart-Smith and Berest-Etingof-Ginzburg (generalizing results for
complex curves and surfaces and rational Cherednik algebras).  We also
use this equivalence to enlarge the category of induced $\D$-modules
on a smooth variety $X$ by collecting induced $\D_X$-modules on
varying cuspidal quotients.  The resulting {\em cusp-induced}
$\D_X$-modules possess both the good properties of induced
$\D$-modules (in particular, a Riemann-Hilbert description) and, when
$X$ is a curve, a simple characterization as the generically
torsion-free $\D_X$-modules.

\end{abstract}

\maketitle

\section{Introduction}
\subsection{$\D$-Modules on Singular Varieties.}
Let $Y$ denote a variety over a field and let $\D_Y$ denote the full
sheaf of differential operators on $Y$ (in characteristic zero $\D_Y$
is the familiar sheaf of differential operators, and in characteristic
$p$ the sheaf $\D_Y$ will contain all divided powers of operators
$\del$).  It is well known that on a {\em general} singular variety
$Y$, the ring of differential operators is badly behaved: in
particular, it need not be Noetherian nor, if $Y$ is affine, a simple
ring (see \cite{BGG}).

Following Grothendieck, there is an alternative notion of a sheaf with
``infinitesimal parallel transport'' called a {\em (co)stratification}
(see \cite{Be1,Be2,Hecke} or Section \ref{Descent and D-modules}),
which agrees with the notion of left (respectively, right)
$\D_Y$-module when $Y$ is smooth but will not, in general, when $Y$ is
singular. Our first result presents a simple vanishing condition that
guarantees that the notions of $\D$-module and (co)stratification on a
singular variety coincide.

\noindent
\begin{thm}\label{D vs strat}  Let $Y$ be a good Cohen-Macaulay variety
(see Definition \ref{good CM}). Then the categories of left
$\D_Y$-modules and stratifications on $Y$ are naturally equivalent, as
are the categories of right $\D_Y$-modules and costratifications on
$Y$.
\end{thm}

The class of good Cohen-Macaulay varieties, characterized by the
vanishing of higher local cohomology sheaves along the diagonal,
includes the smooth varieties and those for which the diagonal is a
set-theoretic local complete intersection, in particular varieties
with cusp singularities. Our technique is based on an extension to
this setting (Theorem \ref{homfromjetsandtensorwithD}) of the
Grothendieck-Sato formula
\begin{equation*}
\D_X=\underline{H}^d_{\Delta}(X\times X,\Oo_X\boxtimes\omega_X)
\end{equation*}
describing the sheaf of differential operators on a smooth variety $X$
of dimension $d$ as a local cohomology sheaf along the diagonal. This
formula is a generalization of the description of differential
operators on a smooth curve as kernel functions with poles along the
diagonal, through the Cauchy integral formula (see e.g. \cite{BS}).
The theorem then follows from the observation that $\D_Y$ is {\em flat}
over $\OY$ whenever the variety $Y$ is a good Cohen-Macaulay variety
in the sense of Definition \ref{good CM}.

\subsection{Cusp Morita Equivalence}
We will refer to a map $X\rightarrow Y$ as a {\em cuspidal quotient}
morphism when $Y$ is the quotient of $X$ by an infinitesimal
equivalence relation (so that $f$ is a homeomorphism on underlying
ringed spaces---we give a precise definition in Section \ref{cusp
quotients and jets}). Examples include the normalization map of a
curve with cusp singularities, the normalization map $\h\to X_m$ of
the space of quasiinvariants for a Coxeter group (see below), and the
Frobenius homeomorphism in characteristic $p$. A good cuspidal
quotient is one which satisfies a relative version of the good
Cohen-Macaulay condition (vanishing of local cohomology along the
graph), and hence benefits from a relative version of the
Grothendieck-Sato formula. By exploiting the resulting identification
of $\D$-modules and stratifications, we obtain the following:

\noindent
\begin{thm}\label{cusp Morita}
Let $f:X\to Y$ be a good cuspidal quotient morphism between good
Cohen-Macaulay varieties (in particular, any cuspidal quotient
morphism from a smooth variety $X$ to a CM variety $Y$).  Then $\D_X$
and $\D_Y$ are canonically Morita equivalent.
\end{thm}

\begin{corollary}
If $X\rightarrow Y$ is a universal homeomorphism from a smooth
 variety $X$ to a CM variety $Y$, then $\D_Y$ is (left and right)
 Noetherian.  If $Y$ is affine, then $\D(Y)$ is a simple ring.
\end{corollary}
\noindent
More precisely, we show that the usual description (in the smooth
setting) of $\D$-module pushforward and pullback by bimodules $\DXY$
and $\DYX$ may be adapted to this singular setting, and that the
functors between left and right $\D$-modules given by these bimodules
are equivalences.

$\D$-modules on singular curves over $\C$ (or an algebraically closed
field of characteristic zero) have been extensively studied
\cite{Smith,Mu,DE,BW ideals}; in particular, Smith and Stafford proved
\cite{SS}, using techniques of noncommutative algebra, that the
category of $\D$-modules on a cuspidal curve is Morita equivalent to
the category of $\D$-modules on its (smooth) normalization. This
result was extended to a class of singular surfaces by Hart and Smith
\cite{HS}, who conjectured an extension to arbitrary Cohen-Macaulay
varieties with a smooth bijective normalization; further refinements
also appeared in \cite{CS,Jones}.  Higher-dimensional examples of
cuspidal quotients are given by the normalization map of the variety
of $m$-quasiinvariants $X_m$ for a Coxeter group $W$ acting by
reflections on a vector space $\h$.  The variety $X_m$ sits in a
diagram $\h\to X_m\to \h/W$ with $\h\to X_m$ bijective, and arises in
the study of the quantum Calogero-Moser dynamical system; it was
proven in \cite{BEG} that $\cD(X_m)$ is simple.  Theorem \ref{cusp
Morita} generalizes and clarifies these results of \cite{SS,HS,BEG}
(and resolves the conjecture of \cite{HS}).  In the case of the
inspiring but atypical {\em flat} example of a cusp quotient, the
Frobenius homeomorphism, our result becomes the Cartier descent for
stratifications due to Berthelot \cite{Be2}.

\subsection{Crystals on Cusps}

The notion of right $\D$-module or costratification may be abstracted
further into the definition of a $!$-{\em crystal} (on the
infinitesimal site---see \cite{Hecke}, 7.10), which is a sheaf endowed
with compatible extensions (``parallel transport'') to arbitrary
nilpotent thickenings. Crystals have many good properties, and in
particular are characterized by the ``Kashiwara theorem'': for any
closed embedding of $Y$ in a smooth variety $Z$, $!$-crystals on $Y$
are just right $\D_Z$-modules supported on $Y$.  Generally,
$\D_Y$-modules and costratifications also diverge from crystals on
$Y$.  However, we explain in Proposition \ref{crystals descend} that
when $X$ is smooth, it is a formal consequence of the Kashiwara
theorem that the categories of $!$-crystals on $X$ and its cuspidal
quotients are equivalent. When combined with the above descent for
$\D$-modules, this implies that right $\D_Y$-modules are equivalent to
$!$-crystals on $Y$ and therefore satisfy the Kashiwara theorem.
\begin{corollary}[Cusp Kashiwara Theorem]
Suppose that $X\rightarrow Y$ is a universal homeomorphism from a smooth
variety $X$ to a CM variety $Y$.  Then right $\D_Y$-modules (or
costratifications) are equivalent to $!$-crystals on $Y$ and thus are
equivalent to right $\D_Z$-modules supported on $Y$ for any closed
embedding $Y\hookrightarrow Z$ of $Y$ in a smooth variety $Z$.
\end{corollary}

 This generalizes the result of \cite{DE} for curves, and resolves
Conjecture 9.9 of \cite{BEG} for arbitrary cusps (we note that this
conjecture in the original case of Cherednik algebras follows from the
Morita equivalence of \cite{BEG} and the standard descent for crystals
of Proposition \ref{crystals descend}).

The intuition behind all these results comes from the description of
stratifications as sheaves equipped with the structure of
infinitesimal parallel transport, that is, equivariance for the deRham
groupoid on $X$ given by the formal neighborhood of the diagonal in
$X\times X$.  This means that the structure of stratification (or
$\D$-module in the smooth case) may be interpreted as a kind of
descent datum, indicating how to descend a sheaf on $X$ to the
quotient of $X$ by the deRham groupoid (known as the deRham space of
$X$).  On the other hand, cuspidal quotients $X\to Y$ are precisely
the quotients of $X$ by {\em subgroupoids} of the deRham groupoid
(that is, by an equivalence relation living along the diagonal in
$X\times X$). In other words, a stratification on $X$ already comes
equipped with descent data for any cuspidal quotient $Y$. Further, one
may hope that the descended sheaf retains enough of the infinitesimal
structure from $X$ to be itself a stratification on $Y$.  The cusp
structures, which are certain slight shrinkings of the sheaf of
functions on $X$, may be imagined to arise by letting the smooth
variety $X$ ``drip'' or ``pinch'' a little; one thereby obtains a
system of ``dripping varieties'' all of which are dripping down toward
the deRham space of $X$, hence have the same collection of
$\D$-modules (which are, morally speaking, sheaves pulled back from
this deRham space).

Unfortunately, the naive $\Oo$-module descent is not the correct
adjoint to the pullback functor in the category of deRham-equivariant
sheaves---in fact, effective descent fails for $\Oo$-modules in the
cusp setting. However, perhaps surprisingly this descent picture (made
precise in a suitable way) {\em does} still apply for the category of
deRham-equivariant sheaves, giving the equivalence of categories of
$\D$-modules of Theorem \ref{cusp Morita}.

\subsection{Cusp Induction.}
The Morita equivalence for cusps has applications to the study of
$\D$-modules on a smooth variety $X$ as well. There is an important
subcategory of the category of (left or right) $\D_X$-modules
consisting of {\em induced} $\D$-modules (see \cite{S,chiral}): an
induced right $\D$-module is one of the form ${\mathcal
F}\otimes_{\theo_X}\D_X$ for a quasicoherent $\theo_X$-module
${\mathcal F}$.  This category generates the derived category of
$\D_X$-modules and is convenient to work with in many ways, but it is
rather small and seems to have no known intrinsic characterization as
a subcategory of the category of $\D_X$-modules.  However, we can use
the Morita equivalence of Theorem \ref{cusp Morita} to ``collect'' the
categories of induced $\D$-modules from all cuspidal quotients of $X$.
It is tempting to think of this construction as a substitute for
Cartier descent in characteristic $p$, defining integrable connections
on $X$ by pullback of quasicoherent sheaves under powers of Frobenius:
in the absence of this canonical cofinal collection of cusp quotient
maps in characteristic zero, we use the collection of all cusps.  By
collecting modules we obtain a larger category sharing all the good
properties of induced $\D$-modules and, in the case of curves,
admitting a simple intrinsic characterization. We call a $\D_X$-module
{\em cusp-induced} if it lies in the essential image of the category
of induced $\D_Y$-modules under the equivalence of Theorem \ref{cusp
Morita} for some Cohen-Macaulay cuspidal quotient $X\rightarrow Y$. We
thereby obtain an equivalence of categories:

\noindent
\begin{thm}\label{induction theorem}
There is an equivalence of categories between the direct limit
of $\OY$-modules with differential operators as
morphisms and cusp-induced $\D$-modules:
$$(\underset{\underset{X\rightarrow Y}{\longrightarrow}}{\lim}
\operatorname{qcoh}(\theo_Y),\on{Diff})\to \on{cusp-ind}(\cD_X).$$ A
quasi-inverse functor is given by the deRham functor.
\end{thm}

This theorem may be considered as a ``cuspidal Riemann-Hilbert
correspondence'', describing the deRham functor on the full
subcategory of cusp-induced $\D$-modules. In the case when $X$ is a
smooth curve, the category of cusp-induced $\D$-modules is precisely
identified with the category of $\D_X$-modules that are generically
torsion-free (see Proposition \ref{cusp-induced on a curve}). In
particular, this applies to (locally) {\em projective} $\D$-modules,
called $\D$-{\em bundles} in \cite{chiral}. As we explain in Sections
\ref{cusp-induced modules section}, \ref{cusp RH section}, and
\ref{cusp-induced on curve section}, the theorems above reproduce and
generalize results of Cannings and Holland \cite{CH ideals, CH cusps}
describing $\D$-bundles on a nonsingular curve. More precisely,
$\D$-bundles are classified by their deRham data, which are
torsion-free sheaves on deep enough cusps (a rigified form of the
deRham datum, known as a ``fat sheaf,'' is used in \cite{BGK2} to
study quiver varieties).  The deRham data are parametrized by an
infinite-dimensional Grassmannian---the {\em adelic Grassmannian}
$\Gr^{ad}(X)$ of \cite{Wilson}---which may be described as the direct
limit of the compactified Picard varieties (moduli of rank 1
torsion-free sheaves) of the dripping curves $Y$. Note that the
projective rank 1 {\em induced} $\D$-modules are classified simply by
$\operatorname{Pic}(X)$; more generally, a submodule of a locally
free, or even cusp-induced, $\D$-module with a finitely supported
cokernel is cusp-induced in any dimension.

The adelic Grassmannian and related moduli spaces of projective
$\D$-modules (``$\D$-bundles,'' \cite{chiral}) have appeared
recently in several contexts \cite{BGK1,BGK2,BW automorphisms, BW
ideals}. In particular, the adelic Grassmannian was introduced by
Wilson (\cite{Wilson} for $X={\mathbb A}^1$) to collect the data for
the algebraic (Krichever) solutions to the KP hierarchy coming from
all cusp quotients $Y$ of $X$, which we see are naturally described by
$\D$-line bundles on the smooth curve $X$. In \cite{solitons}, we give
a completely different (morally ``Fourier dual'') construction of
solutions of the KP hierarchy from $\D$-bundles, in particular
explaining the mysterious link between solitons and many-body
(Calogero-Moser) systems discovered in \cite{AMM,Kr1,Kr2} and deepened in
\cite{Wilson}. In \cite{W}, we use $\D$-bundles (specifically, a
factorization structure on the adelic Grassmannian) to give a
geometric construction of the $\cW_{1+\infty}$-vertex algebra
(associated to the central extension of the Lie algebra of
differential operators on the circle) and localization for its
representations.

The reader may also wish to see \cite{BW diff} for a recent review of
a problem closely related to the subject of the present paper, the
problem of classifying affine varieties up to differential
isomorphism, where cuspidal quotients of affine curves play an
important role.

\subsection{Overview and Acknowledgements.}  In Section \ref{Descent and
D-modules} we prove the Grothendieck-Sato description of differential
operators in the generality we require.  In Section \ref{jets and
stratifications} we use the Grothendieck-Sato formula to compare the
categories of $\D$-modules and stratifications. Section \ref{cusp
morita equivalence section} contains the proof of the Morita
equivalence for $\D$-modules under cuspidal quotients. Finally in
Section \ref{cusp induction section} we introduce cusp-induced
$\D$-modules and describe their main properties.

The authors are grateful to Brian Conrad, Pierre Deligne,
 Victor Ginzburg, Ian
Grojnowski, Ernesto Lupercio, Tony Pantev and especially Matthew
Emerton, Dennis Gaitsgory, Toby Stafford and Michel Van den Bergh for
conversations concerning this work.  Both authors were supported
in part by
NSF postdoctoral fellowships and MSRI postdoctoral fellowships.

\section{Grothendieck-Sato Formula for Differential Operators}\label{Descent 
and D-modules}
\begin{convention}
Fix a ground field $k$ of any characteristic.  For us, all schemes and
morphisms are defined over $k$.  By a {\em variety} we will mean an
integral, separated scheme of finite type over $k$.
\end{convention}

\subsection{Cuspidal Quotients and Jets}\label{cusp quotients and jets}

Recall that a morphism of schemes
 $f:X\rightarrow Y$ is a {\em universal homeomorphism} if
for every morphism $Y'\rightarrow Y$ the pullback
$f_{(Y')}: Y'\times_Y X \rightarrow Y'$ is a homeomorphism.
It follows from Prop. IV.2.4.5 of \cite{EGA} that a morphism
$f:X\rightarrow Y$ of $k$-varieties is a universal homeomorphism
if and only if $f$ satisfies:
\begin{enumerate}
\item $f$ is a finite morphism.
\item $f$ is surjective.
\item For every algebraically closed field $K$, the map
$X(K)\xrightarrow{f(K)} Y(K)$ is injective.
\end{enumerate}

\begin{example}
If $Y$ is a $k$-variety, then the normalization map
$\widetilde{Y}\rightarrow Y$ is a universal homeomorphism if and only if
$Y$ is geometrically unibranch (see Section IV.6.15 of \cite{EGA})---this 
follows from Lemme 0.23.2.2 of \cite
{EGA}.

Note that a variety is geometrically unibranch if and only if
it is \'etale-locally irreducible (this is a consequence of
\cite[IV.18.8.15]{EGA}).
\end{example}

\begin{defn}
A {\em cuspidal quotient morphism} $X\rightarrow Y$ is a universal
homeomorphism between Cohen-Macaulay varieties $X$ and $Y$ over $k$.
\end{defn}

\begin{example}
Suppose $X$ is a nonsingular variety defined over a field $k$ of
characteristic $p$.  Then the Frobenius morphism $X\rightarrow X'$ is
a cuspidal quotient morphism
(see Exp. XV, Section 1, Prop. 2 of \cite{SGA5}).
\end{example}

\begin{defn}[Jets]\label{jet defn}
For a variety $X$ over $k$, write $\cJ_X = \theo_{\widehat{X\times
X}}=\underset{\longleftarrow}{\lim}\,\theo_{X\times X}/I^k_{\Delta}$,
considered as a pro-coherent $\OX$-bimodule.  The bimodule $\cJ_X$,
which is actually a Hopf algebroid, is called the {\em jet algebroid}
of $X$.  The associated formal groupoid $\widehat{\Delta} =
\operatorname{FSpec}(\cJ_X)$ is the {\em jet groupoid} of $X$.  More
generally, for any finite morphism morphism $f:X\rightarrow Y$ we
write $\JYX = \underset{\longleftarrow}{\lim}\,\theo_{Y\times
X}/I^k_{\Gamma_f}$ for the ring of functions on the formal completion
of $Y\times X$ along the graph $\Gamma_f$ of $f$, and $\JXY$ for the
completion of the graph in $X\times Y$.
\end{defn}

\begin{remark}[Characteristic $p$]
As we mentioned in the introduction, in characteristic $p$ we work always
with the full formal neighborhood of the diagonal as our jet groupoid,  with 
the result
that the ring $\D$ of differential operators that we consider is the full ring 
of
differential operators (and in particular includes all divided powers).
This has the advantage that it allows us
to use Kashiwara's theorem in Sections \ref{cusp morita equivalence section} 
and
\ref{cusp induction section} (which would not hold if we used a smaller 
subgroupoid).
\end{remark}

\begin{remark}[Pro-Coherent Sheaves]  The jet algebroid is not
quasicoherent, but rather {\em pro-coherent}, that is, it is a limit
of a filtered system of coherent sheaves; in working with pro-coherent
objects, one remembers the inverse system from which the limit arose,
up to a relaxed notion of isomorphism (see the appendix of \cite{Artin-Mazur} 
or Deligne's appendix to \cite{RD} for
details).  In working with these algebroids, all operations (tensor
product, $\uHom$, etc.) are taken in the category of
pro-quasi-coherent sheaves, where the quasi-coherent sheaves are taken
as constant inverse systems (in topological terms, with the discrete
topology).  So, for example, if $\{\cJ^n\}_{n\geq 0}$ is an inverse
system of coherent $\OX$-modules, the tensor product of
$\{\cJ^n\}_{n\geq 0}$ with the quasicoherent $\OX$-module $M$
corresponds to the inverse system $\{M\ot_\OX \cJ^n\}$.

One may define appropriate notions of Hopf algebroid as well as quasicoherent 
comodule and
(as we shall define later) cocomodule for a Hopf algebroid in the pro-
quasicoherent
category: here the required modification is that the morphisms for
the coaction, counit, etc. are morphisms in the pro-category.
\end{remark}

\subsection{$\cD$-Modules and Local Cohomology.}
\label{section on jets and D-modules}
The Grothendieck-Sato formula describes the relationship
between $\D_Y$ and the formal completion of the diagonal in
$Y\times Y$, as we will explain.

\begin{prop}\label{CMloccohvanishing}
Suppose $X$ is a Cohen-Macaulay variety over $k$ and $f:X\rightarrow
Y$ is a finite morphism.  Let $d = \dim(Y) = \dim(X)$.  For any
quasicoherent $\OY$-module $M$, one has \bd \uExt_Y^i(\JYX, M) =
\underline{H}^{i+d}_{\Gamma_f}(Y\times X, M\boxtimes \omega_X) \ed for
all $i\in{\mathbf Z}$.  In particular, the local cohomology sheaf
$\underline{H}_{\Gamma_f}^{d+i}(Y\times X, M\boxtimes \omega_X)$
vanishes for all $i <0$.
\end{prop}
\begin{proof}
The formula follows from Grothendieck duality.  By assumption, $X$ is
of finite type; hence by Nagata's Theorem there is an embedding
$X\overset{i}{\hookrightarrow}\overline{X}$ in a proper $k$-scheme
$\overline{X}$.  Let $\Gamma_f^{(k)}$ denote the $k$th-order
neighborhood of the graph of $f$ in $Y\times\overline{X}$.  Let $\pi:
Y\times\overline{X}\rightarrow Y$ denote projection on the first
factor.  Then for any quasicoherent $\OY$-module $M$, Grothendieck
duality implies that \bd
\bR\pi_*\left(\bR\uHom_{Y\times\overline{X}}\left(\theo_{\Gamma_f^{(k)}},
\pi^!M\right)\right) \simeq
\bR\uHom_Y\left(\bR\pi_*\theo_{\Gamma_f^{(k)}}, M\right).  \ed Because
$\theo_{\Gamma_f^{(k)}}$ is supported on the subscheme
$\Gamma_f^{(k)}\subset Y\times\overline{X}$ that is finite over $Y$
and on the scheme $Y\times X$ which is Cohen-Macaulay over $Y$, this
reduces (see \cite{Conrad}) to \bd \pi_*\uExt^{i+d}_{Y\times
X}(\theo_{\Gamma_f^{(k)}}, M\boxtimes \omega_X) \simeq
\uExt_Y^i(\pi_*\theo_{\Gamma_f^{(k)}}, M).  \ed Taking colimits over
$k$, the left-hand side becomes local cohomology, giving the desired
formula and vanishing.
\end{proof}

We next wish to determine when the local cohomology groups of
Proposition \ref{CMloccohvanishing} vanish for $i>0$; unfortunately,
this is a slightly more complicated question.
\begin{defn}[Good morphism, good CM variety] \label{good CM}
A {\em good} cuspidal quotient morphism
 $f:X\rightarrow Y$ between
Cohen-Macaulay varieties of dimension $d$ is a cuspidal quotient
morphism with the property that for any quasicoherent $\OY$-module $M$, the
local cohomology sheaf $\underline{H}_{\Gamma_f}^{d+1}(Y\times X,
M\boxtimes\omega_X )$ vanishes.  The cuspidal quotient morphism $f$ is
{\em very good} if for any quasicoherent $\theo_{Y\times X}$-module $M$,
the local cohomology sheaf $\underline{H}^{d+i}_{\Gamma_f}(Y\times X,
M)$ vanishes for all $i>0$.  A CM variety $Y$ of dimension $d$ is {\em
good} ({\em very good}) if the identity morphism $\operatorname{id}_Y:
Y\rightarrow Y$ is a good (respectively, very good) morphism.
\end{defn}

\noindent
Of course, a very good cuspidal quotient morphism is good.
\begin{defn}
We will call an integral closed subscheme $Z\subset W$ of an integral
scheme $W$ a {\em set-theoretic local complete intersection} if, for
each point $z\in Z\subset W$, there is a regular sequence $x_1, \dots,
x_k$ in $\theo_{W,z}$ such that $\spec\theo_{W,z}/(x_1,\dots, x_k)$ is
a nilpotent thickening of its closed subscheme
$\spec\theo_{Z,z}$.
\end{defn}
We next give some conditions that imply that a variety or morphism
is very good.
\begin{prop}\label{loccohvanishingoncusp}\label{CMdiagislci}
\mbox{}
\begin{enumerate}
\item If $X$ is a smooth $k$-variety of dimension $d$, then $X$
is very good.
\item Suppose that $f:X\rightarrow Y$ is a cuspidal quotient morphism.
Then the following are equivalent:
\begin{enumerate}
\item $X$ is very good.
\item $f$ is very good.
\item $Y$ is very good.
\end{enumerate}
\item Suppose that $f:X\rightarrow Y$ is a cuspidal quotient morphism
 and that
 $\Gamma_f\subset Y\times X$ is a set-theoretic local complete
intersection in $Y\times X$.  Then $f$ is very good.
\end{enumerate}
\end{prop}
\begin{proof}
Part (1) is well-known: it follows from the fact that $\Delta_X$ is a
local complete intersection in $X\times X$ when $X$ is smooth over
$k$.  For part (2), because we wish to compute local cohomology sheaves along 
$\Delta_X$, $\Delta_Y$ and
the graph of $f$, we may assume that $X$ and $Y$ are affine schemes.  Let
$U_Y\subset Y\times Y$ (respectively $U_X$, respectively $U_f$) denote the 
complement of
$\Delta_Y$ in $Y\times Y$ (respectively of $\Delta_X$ in $X\times X$, 
respectively of $\Gamma_f$ in $Y\times X$).
It follows from the assumption that $X\rightarrow Y$ is a cuspidal
quotient morphism that the natural maps
$X\times X\xrightarrow{f\times 1} Y\times X \xrightarrow{1\times f}
 Y\times Y$ are homeomorphisms and
that the maps $\Delta_X\rightarrow
X\times_Y X\rightarrow \Delta_Y$ are bijective, so $U_X = (f\times 1)^{-1}
(U_f) =  (f\times f)^{-1}(U_Y)$ and $(1\times f)^{-1}(U_Y) = U_f$; moreover 
the maps
$U_X\rightarrow U_f\rightarrow U_Y$ are surjective and finite.
By Proposition 1.1 of
\cite{HartshorneCD}, this implies that the cohomological dimensions of
$U_X$, $U_f$ and $U_Y$ are the same.
By Proposition 2.2 of \cite{HartshorneLC}, one has $H^i(U,M) =
H^{i+1}_{\Gamma}(M)$ for every quasicoherent $M$ and every $i\geq 1$
(here $\Gamma$ stands for one of $\Delta_X$, $\Delta_Y$ or $\Gamma_f$
and $U$ stands for the corresponding open complement).  Since we have
reduced to the case of an affine ambient variety, the local cohomology
sheaf is the sheaf associated to the module $H^{i+1}_{\Gamma}(M)$, and
consequently the maximal degrees for nonzero local cohomology sheaves
along $\Delta_X$, $\Delta_Y$ and $\Gamma_f$ all coincide.

For part (3),
the question is local on $Y$, so we may assume that there is a regular
sequence $x_1,\dots, x_d$ in $\theo_{Y\times X}$ such that
$\spec\theo_{Y\times X}/(x_1,\dots,x_d)$ is a nilpotent thickening of
$\Gamma_f$.  But then the local cohomology along $\Gamma_f$ coincides
with local cohomology with respect to $(x_1,\dots, x_d)$, which is
computed by the direct limit of Koszul cohomology groups by Theorem
2.3 of \cite{HartshorneLC}, hence vanishes above degree $d$.
\end{proof}

Recall (following Grothendieck) that
$\cD_Y = \uHom_\OY(\cJ_Y,\OY)$, and more generally
$\cD_Y(M,N) = \uHom_\OY(\cJ_Y\otimes_\OY M, N)$ for any quasicoherent
$\OY$-modules $M$ and $N$; we will
require also the relative version for a map $X\rightarrow Y$,
replacing $\cJ$ by its analogs from Definition \ref{jet defn}:
\begin{defn}\label{def of bimodules}\mbox{}
\begin{enumerate}
\item Let $\DXY=\uHom_\OX (\JXY, \OX)$.
\item Let $\DYX = \uHom_\OY(\JYX, \OY)$.
\end{enumerate}
\end{defn}

\begin{lemma}\label{jetsgivejets}
Let $f:X\rightarrow Y$ be an affine morphism between $k$-varieties.
\begin{enumerate}
\item The pro-coherent sheaf $\JYX$ (respectively $\JXY$)
 associated to the completion of $Y\times X$ ($X\times Y$)
along the graph of $f$ is $(1\times f)^*\cJ_Y$
(respectively $(f\times 1)^*\cJ_Y$).
\item
If $X\xrightarrow{f} Y$ is a universal homeomorphism, then the
jet algebroid of $X$ is the pullback of the jet algebroid of $Y$; that
is, $\cJ_X=(f\times f)^*\cJ_Y$ as pro-coherent Hopf algebroids.
\end{enumerate}
\end{lemma}
\begin{proof}
In (1), the statements for $\JYX$ and $\JXY$ hold because the ideal
of $\Gamma_f$ is generated by elements $a\otimes 1 - 1\otimes a$ where $a$ is 
a local
section of $\OY$.

For (2),
we have $\cJ_X = \underset{\longleftarrow}{\lim}\,\theo_{X\times
X}/I_{\Delta}^k$ and $\cJ_Y =
\underset{\longleftarrow}{\lim}\,\theo_{Y\times Y}/I_{\Delta}^k$.
Locally on either $X$ or $Y$, $I_{\Delta}$ is generated by elements of
the form $a\otimes 1 - 1\otimes a$ for local sections $a\in\theo$.  It
follows that one has an exact sequence \bd (f\times
f)^*I_{\Delta_Y}\rightarrow \theo_{X\times X}\rightarrow
\theo_{X\times_Y X} \rightarrow 0 \ed on $X\times X$.  In particular,
we find that the image of $(f\times f)^*I_{\Delta_Y}$ in
$\theo_{X\times X}$ is contained in $I_{\Delta_X}$. Moreover, since
$f$ is a universal homeomorphism,
 $\theo_{X\times_Y X}$ is supported in a finite-order
neighborhood of the diagonal in $X\times X$.  It thus follows that
the image of $(f\times f)^*I_{\Delta_Y}$
 contains $I_{\Delta_X}^n$
for some $n\geq 1$.  Therefore, the two limits coincide.  Since the
two inverse systems are cofinal, the pro-structures also agree.
\end{proof}

\begin{corollary}\label{sheaves of diff ops}
Suppose that $f:X\rightarrow Y$ is a cuspidal quotient morphism.  Then
$\DXY$ is canonically identified with the sheaf $\D_Y(\OY,\OX)$ of
$\OY$-differential operators from $\OY$ to $\OX$, and similarly
$\DYX=\D_Y(\OX,\OY)$.
\end{corollary}
\begin{proof}
By part (1) of Lemma \ref{jetsgivejets} we have $\JXY = \OX\otimes\cJ_Y$, so
we get $\DXY = \uHom_{\OX}(\OX\otimes \cJ_Y,\OX)$; the right-hand side equals
$\uHom_{\OY}(\cJ_Y,\OX)$ by the usual adjoint associativity, and this is
$\cD_Y(\OY,\OX)$ by Grothendieck's definition of differential operators.

Part (1) of Lemma \ref{jetsgivejets} also gives $\JYX = \cJ_Y\otimes \OX$, so 
the
equation $\DYX = \cD_Y(\OX,\OY)$ follows from Grothendieck's definition.
\end{proof}
See Lemma \ref{pullbacks of D} and Remark \ref{the D bimodules} for an
alternative definition of $\DXY$ and $\DYX$ and their relation with
the standard bimodules of $\D$-module theory.

\begin{thm}[Grothendieck-Sato Formula]\label{homfromjetsandtensorwithD}
Suppose $f:X\rightarrow Y$ is a good cuspidal quotient morphism of
$k$-varieties.  Let $d=\dim(Y) = \dim(X)$.  Then
for any quasicoherent $\OY$-module $M$,
\begin{equation}\label{GroSat}
\underline{H}_{\Gamma_f}^d(Y\times X, M\boxtimes\omega_X) =
\uHom_{\OY}(\JYX, M) = M\ot_{\OY}\DYX.
\end{equation}
\end{thm}

\begin{proof}
The first equality of \eqref{GroSat} is immediate from Proposition
\ref{CMloccohvanishing}.  To prove the second equality of
\eqref{GroSat}, observe that the
vanishing of $\underline{\operatorname{Ext}}^1_{\OY}(\JYX, M) =
\underline{H}_{\Gamma_f}^{d+1}(Y\times X, M\boxtimes\omega_X)$ for
every quasicoherent $M$ implies that $\uHom_{\OY}(\JYX, \,\cdot\,)$ is
exact, and therefore this exactness follows from $f$ being good.
Given a finite presentation $\theo_Y^A\xrightarrow{\alpha}
\theo_Y^B\xrightarrow{\beta} M \rightarrow 0$ of a coherent
$\OY$-module $M$, we get a commutative diagram \bd
\xymatrix{\theo_Y^A\otY\DYX \ar[r]\ar[d]^{\cong} & \theo_Y^B \otY\DYX
\ar[r]\ar[d]^{\cong}& M\otY\DYX \ar[r]\ar[d] & 0\\ \uHom_{Y}(\JYX,
\theo_Y^A)\ar[r] & \uHom_{Y}(\JYX, \theo_Y^B)\ar[r] & \uHom_{Y}(\JYX,
M)\ar[r] & 0} \ed with exact rows, which proves \eqref{GroSat} for coherent
$M$.  Finally, local cohomology \cite[Proposition 1.12]{HartshorneLC} and
tensor product both commute with colimits, implying \eqref{GroSat} for
all quasicoherent $M$.
\end{proof}

\begin{corollary}\label{flatness cor}
For $f:X\to Y$ as in Theorem \ref{homfromjetsandtensorwithD},
the sheaf $\DYX$ is flat over $\OY$. In
particular if $Y$ is a good Cohen-Macaulay variety, then $\D_Y$ is
flat as a left $\OY$-module.
\end{corollary}

Applying the theorem to the identity morphism of a good Cohen--Macaulay 
variety, we
obtain the corollary:

\begin{corollary} The Grothendieck--Sato formula
$$\D_Y=\underline{H}^d_{\Delta}(Y\times Y,\Oo_Y\boxtimes\omega_Y)$$
holds for good Cohen--Macaulay varieties $Y$.
\end{corollary}

\section{Differential Operators, Jets and Stratifications}\label{jets and 
stratifications}

In this section, after reviewing some of the formalism of descent, we
compare left and right $\D$-modules with stratifications and
costratifications, or equivalently jet comodules and jet cocomodules,
on a good Cohen-Macaulay variety $Y$. These comparisons follow from
the extension of the Grothendieck-Sato formula, Theorem
\ref{homfromjetsandtensorwithD} (applied to the identity map of
$Y$). We also explain the descent for crystals under cuspidal
quotients. An excellent reference for stratifications and
costratifications in the smooth setting is provided by \cite{Be1,Be2}.

\subsection{$*$-Equivariance}
For background
on the relationship between groupoids and Hopf algebroids, see, for
example, \cite{Hovey}.

\begin{defn}
Suppose $\cG\rightrightarrows X$ is a groupoid over $X$ (we suppress
the unit in this notation).  We say an $\OX$-module $M$ is {\em
equivariant with respect to the groupoid $\cG\rightrightarrows X$} if
$M$ is equipped with an isomorphism $p_1^*M\to p_0^*M$ that is
compatible with composition over $\cG\times_X\cG$ in the usual sense
(see Section 1.6 of \cite{Deligne}).
\end{defn}
We work with groupoids that are affine or pro-affine over $X$, and
write $\cH$ for the Hopf algebroid associated to $\cG$.  By the
adjunction of pushforward and pullback, we have
$$\uHom_{\theo_\cG}(p_1^*M, p_0^*M)= \uHom_{\Oo_X}(M,p_{1*}p_0^*M).$$
We identify $p_{1*}p_0^*M=M\otX \cH$, considered as an $\OX$-module
via the right $\OX$-structure on $\cH$ and where the tensor product is
taken over the left structure.
\begin{lemma}\label{comod}
The category of $\cG$-equivariant modules is equivalent to the
category of right $\cH$-comodules that are counital in the sense that
$$M\to M\otX \cH\xrightarrow{\on{counit}}M\otX \OX=M$$ is the
identity.
\end{lemma}
\begin{remark}
The counital condition on $\cH$-comodules is equivalent to the requirement that
the corresponding map $p_0^*M\to p_1^*M$ is an isomorphism.
\end{remark}

\begin{prop}\label{pullbackofHopfalgebroids} Suppose $f:X\rightarrow Y$ is 
affine.  Let ${\mc H}_Y$ be a Hopf algebroid on $Y$,
and ${\mc H}_X=(f\times f)^*{\mc H}_Y$.  Then ${\mc H}_X$ is a Hopf
algebroid on $X$ and $f^*$ induces a functor from $\on{comod}({\mc
H}_Y)$ to $\on{comod}({\mc H}_X)$.
\end{prop}

\begin{proof}
The Hopf algebroid structure on ${\mc H}_X$ is defined as follows: the 
comultiplication on ${\mc H}_X=\OX\otY{\mc H}_Y\otY\OX$ is given by
means of the comultiplication on ${\mc H}_Y$ and the map
${\mc H}_Y\otY{\mc H}_Y\to {\mc H}_Y\otY\OX\otY{\mc H}_Y$ induced from the
structure morphism $\OY\to\OX$. The other structures of Hopf algebroid
are immediate, as is the functor on comodules.
\end{proof}

\subsection{!-Equivariance}
We now replace the adjoint pair of functors $(f^*,f_*)$ by the adjoint
pair $(f_*,f^!)$.
\begin{defn}[See Section 7.10 of \cite{Hecke}]
For a finite
morphism of schemes $f$, or more generally for an ind-finite morphism
of formal schemes, the pushforward $f_*:\OX\mbox{-}\on{mod}\to\OY\mbox{-}\on
{mod}$
has a right adjoint $f^!$, defined as follows:  for an $\OY$-module $N$, $f^!
N$ is
the $\OX$-module corresponding to the $f_*\OX$-module
$\uHom_{\OY}(f_*\OX,N)$ on $Y$.
\end{defn}
 In our situation, $f$ is affine, so we will omit the notation $f_*$ for the 
pushforward unless we wish to emphasize that
we are forgetting the $\OX$-module structure down to $\OY$.

Let $N$ be an $\Oo_X$-module which is $!$-equivariant with respect to
the groupoid $\cG$ over $X$: in other words, we are given an
isomorphism $p_0^! N\to p_1^! N$ compatible with compositions.  By the
adjunction of $!$-pullback and pushforward, we have
\bd
\uHom_{\cG}(p_0^!N, p_1^!N)= \uHom_{\Oo_X}(p_{1*}p_0^!N,N).
\ed
It follows that the analog of Lemma \ref{comod} identifies the
notion of $!$-equivariant sheaf with that of {\it cocomodule}:

\begin{defn}[Cocomodules]
Let $\cH$ be a Hopf algebroid on $X$.  An $\cH$-cocomodule $N$ is an
$\Oo_X$-module $N$, equipped with a morphism $\uHom_{\Oo_X}(\cH,N)\to
N$, so that the two compositions $\uHom_{\Oo_X}(\cH\ot\cH ,N)\to N$
given by the coproduct
$$\uHom_{\Oo_X}(\cH\ot\cH ,N)\stackrel{\Delta^*}{\longrightarrow}\uHom_{\Oo_X}
(\cH,N)\to N$$
and the induced map
$$\uHom_{\Oo_X}(\cH\ot\cH ,N)=\uHom_{\Oo_X}(\cH,\uHom_{\Oo_X}(\cH,N))\to
\uHom_{\Oo_X}(\cH,N)\to N$$ agree.
We further require $N$ to be unital, so that the composition
$$N=\uHomX(\OX,N)\xrightarrow{\on{counit}}\uHomX(\cH,N)\to
N$$ is the identity.
\end{defn}

\begin{remark} If $\cH$ is an $\OX$-coalgebra that is {\em projective} over 
$\OX$,
with dual algebra $\cH^*=\uHom_\OX(\cH,\OX)$, then
$\uHom_{\Oo_X}(\cH,M)=M\ot \cH^*$ and the structure of
$\cH$-cocomodule is equivalent to that of right $\cH^*$-module.
\end{remark}

\begin{prop} Suppose $f:X\rightarrow Y$ is finite.
 Let ${\mc H}_Y$ be a Hopf algebroid on $Y$ and ${\mc H}_X$ the corresponding
algebroid on $X$, ${\mc H}_X=(f\times f)^*{\mc H}_Y$. Then $f^!$
induces a functor from ${\mc H}_Y$-cocomodules to ${\mc
H}_X$-cocomodules.
\end{prop}

\begin{proof}
Let $(p_{0,Y},p_{1,Y}):{\mc G}_Y\to Y\times Y$ be the groupoid
corresponding to ${\mc H}_Y$ and $(p_{0,X},p_{1,X}):{\mc G}_X\to X\times
X$ be that corresponding to ${\mc H}_X$. Thus $f\circ
p_{i,X}=p_{i,Y}\circ (f\times f)$. If $N$ is an ${\mc H}_Y$-cocomodule, that 
is an
$\OY$-module with an isomorphism $p_{0,Y}^!N\to p_{1,Y}^! N$
satisfying a composition rule, we obtain an isomorphism
$$p_{0,X}^!f^!N=(f\times f)^!p_{0,Y}^!N\to (f\times
f)^!p_{1,Y}^!N=p_{1,X}^!f^!N,$$ which defines the desired ${\mc H}_X$-
cocomodule
structure on $f^!N$.
\end{proof}

\subsection{Stratifications and Costratifications Coincide With $\cD$-Modules}

\begin{defn}\mbox{}
\begin{enumerate}
\item A stratification $M$ on $X$ is a quasicoherent sheaf with a
right comodule structure for the jet algebroid on $X$. Equivalently,
$M$ is a sheaf equivariant with respect to the deRham groupoid
$\widehat{X\times X}$, or a sheaf equipped with an
isomorphism $p_1^*M\to p_0^*M$ on $\widehat{X\times X}$ compatible
with composition on the triple product.
\item A costratification $M$ on $X$ is a quasicoherent sheaf with a
cocomodule structure for the jet algebroid on $X$. Equivalently,
$M$ is a sheaf $!$-equivariant with respect to the deRham groupoid
$\widehat{X\times X}$, or a sheaf equipped with an
isomorphism $p_0^!M\to p_1^!M$ on $\widehat{X\times X}$ compatible
with composition on the triple product.
\end{enumerate}
\end{defn}

\begin{remark}[Left and Right]
It is clear from the description in terms of equivariance for a groupoid that 
the datum of a counital right
comodule is equivalent to the datum of a counital left comodule.  When one 
dualizes a comodule structure, however,
to obtain a module over the algebra dual to the Hopf algebroid, one makes a 
choice---in our case, the algebra $\cD$
is the {\em left} dual of jets---which breaks the symmetry.  This explains our 
insistence on the use of right comodules.
\end{remark}

\begin{thm}\label{mainequivthm}
Suppose $Y$ is a \good Cohen-Macaulay variety.
\begin{enumerate}
\item\label{coco equals right D} The categories of quasicoherent
$\cJ_Y$-cocomodules (or costratifications on $Y$) and right
$\D_Y$-modules are equivalent.
\item The categories of quasicoherent right $\cJ_Y$-comodules (or 
stratifications on
$Y$) and left $\cD_Y$-modules are equivalent.
\end{enumerate}
\end{thm}
\noindent
Recall that quasicoherence for $\cJ_Y$-cocomodules and $\D_Y$-modules
means quasicoherence over $\OY$.

\begin{proof}
The first assertion, the equivalence of $\cJ_Y$-cocomodules and right
$\D_Y$-modules, is immediate from the
formula
$\uHom_{\OY}(\cJ_Y, M) = M\ot_{\OY}\D_Y$ of \eqref{GroSat} using
standard methods (see \cite[Prop. 1.1.4]{Be2} for the smooth case).

The proof of the second assertion also follows the standard outline
(see \cite[II, 4.1--4.2]{Be cristalline}), using
 a general algebraic fact
that we now explain.
To begin, let $L$ and $R$ be commutative rings, and let
 $Q=\underset{\longleftarrow}{\lim}\, Q_n$ denote a pro-object in the
category of $(L,R)$-bimodules that are finitely generated over $L$.
Let $P=\underset{\longrightarrow}{\lim}\,
\Hom_L(Q_n,L)$ denote the $L$-dual of $Q$.  For simplicity, we assume
that the maps $Q_{n+1}\rightarrow Q_n$ are surjections.
\begin{lemma}\label{duality lemma}
If the natural map $N\otimes_L P\rightarrow \Hom_L(Q,N)$ is an isomorphism
for all $L$-modules $N$,
then the natural map
\begin{equation}\label{flatness eq}
{\mathbf F}: \Hom_L(P\ot_R M, N) \to \Hom_R(M, N\ot_L Q)
\end{equation}
is an isomorphism of $(L,R)$-bimodules
 for all $L$-modules $N$ and $R$-modules $M$.
\end{lemma}
The method of proof is standard (with the usual modifications for the
category of pro-objects):
\begin{proof}[Proof of Lemma]
The hypotheses of the lemma give a canonical element
$1_Q\in \Hom_L(Q,Q) = Q\otimes_L P$.
Using the action of $R$ on $P$, this element
determines a
 map $\delta: R\rightarrow Q\otimes_L P$.
We also have the natural contraction maps
$\on{tr}_n: P_n\ot Q_n \rightarrow L$ for all $n$; by abuse
of notation, we write this system of homomorphisms as
$\on{tr}: P\ot Q \rightarrow L$.
The hypotheses of the lemma give identities
\begin{equation}\label{pro identities}
(1_Q\ot\on{tr})\circ (\delta\ot 1_Q) = 1_Q \;\;\text{and}\;\;
(\on{tr}\ot 1_P)\circ (1_P\ot\delta) = 1_P
\end{equation}
(where the first is an identity in the pro-category).

Recall that the map ${\mathbf F}$ of \eqref{flatness eq} is given by taking
$\phi: P\ot M\rightarrow N$ to the composite
$M\xrightarrow{M\ot \delta} M\ot P\ot Q\xrightarrow{\phi\ot Q} N\ot Q$.
Similarly, one may define a map ${\mathbf F}^{-1}$
by taking $\psi: M\rightarrow N\ot Q$ to the composite
$M\ot P\xrightarrow{\psi\otimes P} N\ot Q\ot P\xrightarrow{N\ot \on{tr}}
N$.
It is then straightforward to check that
${\mathbf F}$ and ${\mathbf F}^{-1}$ are mutually inverse using
\eqref{pro identities}.
\end{proof}
Part (2) of Theorem \ref{mainequivthm} now follows by taking $L=R=\OY$,
 $Q = \cJ_Y$,
$P = \D_Y$, and $M=N$: the hypothesis of Lemma \ref{duality lemma} follows
from Equation \eqref{GroSat}, and so by the lemma we obtain from
\eqref{flatness eq} an isomorphism
\bd
\uHom_{\OY}(\D_Y\ot_{\OY} M, M) = \uHom_{\OY}(M, M\ot_{\OY}\cJ_Y).
\ed
The compatibilities required for module and comodule structures are
immediate from the constructions.\end{proof}

\begin{corollary}
Suppose $f:X\rightarrow Y$ is a cuspidal quotient morphism between
good CM varieties.  Then the naive functors $f^*$ and $f^!$ on
$\OY$-modules determine pullback functors from the categories of left
and right $\D_Y$-modules to the categories of left and right
(respectively) $\D_X$-modules.
\end{corollary}
\noindent
This follows from the existence of the corresponding pullbacks for
jet comodules and cocomodules.

\subsection{Descent for Crystals.}
The theory of $!$-crystals is developed in \cite{Hecke}, section 7.10,
whose definitions we follow.  It is important to note that even in
characteristic $p$ we allow {\em all} nilpotent thickenings (not just
those with fixed divided power structure) in the definition of a
$!$-crystal, and thus our notion does not coincide with the more usual
crystalline terminology in characteristic $p$.

Let $f:X\to Y$ denote a cuspidal quotient of a smooth variety $X$,
$i:Y\to Z$ a closed embedding into a smooth variety $Z$, and $g=i\circ
f:X\to Z$. In this section we describe descent for $!$-crystals,
namely the statement that the functor $f^!$ induces an equivalence of
categories between $!$-crystals on $X$ and $!$-crystals on $Y$. In
fact, $!$-crystals on $X$ are equivalent to right $\D_X$-modules
while, thanks to the Kashiwara theorem for $!$-crystals,
$!$-crystals on $Y$ are equivalent to right $\D_Z$-modules supported
on $Y$. Thus it suffices to prove that the standard $\D$-module
functors define an equivalence of these two categories.
The result seems to be well known, but since we could not find a
reference we sketch a proof for the benefit of the reader. We are
 grateful to Dennis Gaitsgory for very valuable discussions of
$!$-crystals and their properties.

\begin{prop}\label{crystals descend}
The functors $g_*$ and $g^!$ define quasiinverse equivalences of the
categories $\rmod - \D_X$ of right $\D_X$-modules  and $(\rmod - \D_Z)_Y$
of right $\D_Z$-modules  supported on $Y$.
\end{prop}

\begin{proof}
We use the adjunction between $g^!$ and $g^*$ (note that $g$ is
proper). We first prove that the adjunction $Id\to g^! g_*$ is an
isomorphism.  Let $p_1,p_2: X\times_Z X$ denote the two projections
from the descent groupoid of $X$ over $Z$ to $X$. Then we have a
natural base change equivalence $g^!g_*=(p_2)_* (p_1)^!$.  Now by the
Kashiwara equivalence for right $\D_{X\times X}$-modules supported on
the diagonal $\Delta:X\to X\times X$ with right $\D_X$-modules, we
have an isomorphism $(p_1)^!=\Delta_*$. It follows that
$g^!g_*=(p_2)_* (p_1)^!=(p_2)_* \Delta_* = Id$ as desired.

For the converse, the isomorphism property of the adjunction
$g_*g^!\to Id$, we use a flattening stratification of $f:X\to
Y$. Note that for a {\em flat} cuspidal quotient, the equivalence of
categories is a consequence of flat $*$- and $!$-descent for
coherent sheaves---the proof is as in Section 2.4 of \cite{Be2}.  The 
conclusion now follows
by induction: assuming (by way of inductive hypothesis) the isomorphism for
the subcategory of modules
supported on a closed subvariety $V$ which is a union of strata, we wish to add
a stratum $C$ to obtain a closed subvariety $V\cup C$ and check the 
isomorphism there.
But the canonical
(Cousin) decompositions of $g_*g^!M$ and $M$ agree both on $V$ and (by the
stratum-by-stratum equivalence) on $C$; consequently they agree
on all of $V\cup C$.  It
follows that the adjunction morphism is an isomorphism as claimed.
\end{proof}

\section{Cusp Morita Equivalence}\label{cusp morita equivalence section}
In this section, we will prove that tensor products with the standard
bimodules $\DXY$ and $\DYX$ of $\D$-module theory give Morita
equivalences of categories of $\cD$-modules under cuspidal quotient
morphisms.  The functors from $\D_Y$-modules to $\D_X$-modules are
given by $*$ and $!$-pullback, and one should think of the inverse
functors as descent functors in the categories of jet
(co)comodules.

 Throughout this section, we fix a cuspidal quotient morphism
$X\xrightarrow{f} Y$.

\subsection{Morita Equivalence}

\begin{lemma}\label{lemma from mel}
Suppose $f:X\rightarrow Y$ is a finite, dominant map, and that $M$ is a 
nonzero quasicoherent $\OY$-module.  Then $f^!M$ is
nonzero.
\end{lemma}
\begin{proof}
Suppose $M\neq 0$.  By left exactness of $\uHom$, we may assume that
$M=\OY/P$ for some prime $P$ of $\OY$, and may further localize $\OY$
at $P$ (which we suppress in our notation).  Then $f^!M$ becomes
$\uHom_Y(\OX, \OY/P) = \uHom_{\OY/P}(\OX\ot \OY/P, \OY/P)$, the dual
of a vector space over the field $\OY/P$.  Since $\OX$ is supported at
every point of $Y$, $\OX\ot\OY/P$ is nonzero, completing the proof.
\end{proof}

\begin{lemma}\label{pullbacks of D}
Suppose that $f:X\rightarrow Y$ is a cuspidal quotient morphism.
  Then:
\begin{enumerate}
\item $(1\times
f)^!\D_Y = \DYX$.
\item
If $Y$ is good then
 $(f\times 1)^*\D_Y =
\DXY$.
\item If $f$ is good then $(f\times 1)^*\DYX = \D_X$.
\item
$(1\times f)^! \DXY=\D_X$.
\end{enumerate}
Furthermore, if $Y$ and $f$ are good,
then for any quasicoherent $\OY$-module $M$ we have
\begin{equation}\label{shriek eq}
(1\times
f)^!(M\otY\D_Y) = M\otY \DYX.
\end{equation}
\end{lemma}

\begin{proof}
For (1), we have
\begin{align*}
(1\times f)^!\D_Y & =\uHomY(\OX,\D_Y)=\uHomY(\OX,\uHomY(\cJ_Y,\OY))\\
 & =\uHomY(\cJ_Y\otY\OX,\OY)=\DYX
\end{align*}
using the definitions, adjunction, and part (1) of Lemma \ref{jetsgivejets}.

Parts (2) and (3) follow from the identifications
\begin{align*}
(f\times 1)^*\D_Y = \OX\otY\D_Y = \uHomY(\cJ_Y,\OX) = \uHomX(\JXY, \OX),
\end{align*}
\begin{align*}
(f\times 1)^*\DYX = \OX\otY\DYX = \uHomY(\JYX,\OX) = \uHomX(\cJ_X,\OX)
\end{align*}
using the Grothendieck-Sato formula, adjunction and Lemma
\ref{jetsgivejets}.

For part (4), we have
\bd
(1\times f)^!\DXY = \uHom_{\OY}(\OX,\uHom_{\OX}(\JXY,\OX))
 =\uHom_{\OX}(\cJ_X,\OX)=\D_X
\ed
using Lemma \ref{jetsgivejets}.

Finally, \eqref{shriek eq} follows by computing
\begin{align*}
(1\times f)^!(M\otY\D_Y) & =\uHom_{\OY}(\OX,\uHom_\OY
(\cJ_Y,M))=\uHomY(\cJ_Y\otY \OX,M)\\ &=\uHomY(\JYX,M) = M\ot_{\OY}\DYX
\end{align*}
using the Grothendieck-Sato formula.
\end{proof}

\begin{thm}[Cusp Morita Equivalence]\label{the Morita theorem}
Suppose $f:X\rightarrow Y$ is a good cuspidal quotient morphism of
good CM varieties.  Then the bimodules $\DXY$ and $\DYX$ induce Morita
equivalences of the categories of (left or right) $\D_X$-modules and
$\D_Y$-modules.
\end{thm}
\begin{proof}
We will show that $\DXY\ot_{\D_Y} \DYX = \D_X$ and $\DYX\ot_{\D_X}\DXY
= \D_Y$; the equivalences are then given by tensoring with one or the
other of these bimodules.

It follows immediately from the Grothendieck-Sato formula and Lemma
\ref{pullbacks of D} that $\DXY\ot_{\D_Y}\DYX = \D_X$: indeed, we have
\bd \DXY\ot_{\D_Y}\DYX = \OX\otY \D_Y\ot_{\D_Y}\DYX = \OX\otY\DYX =
\D_X.  \ed
Observe that this morphism is just the one obtained by viewing $\DYX$ and
$\DXY$ as sheaves of differential operators (as in
Lemma \ref{sheaves of diff ops}) and taking $\phi\otimes\psi$ to
the composite differential operator $\phi\circ\psi$ from $\OX$ to $\OX$.

Consider the natural map $\DYX\ot_{\D_X}\DXY \xrightarrow{\Phi} \D_Y$,
which again comes from viewing $\DYX$ and $\DXY$ as sheaves of
differential operators and taking composites.  Tensoring on the right
by $\ot_{\D_Y}\DYX$ gives \bd \DYX\ot_{\D_X}\DXY\ot_{\D_Y}\DYX = \DYX
\rightarrow \DYX \ed which is the identity map; letting $M =
\coker(\Phi)$, we thus have $M\ot_{\D_Y}\DYX = 0$.

There is a coequalizer diagram
\begin{equation}\label{coequalizer defining tensor product}
M\otY \D_Y\otY \D_Y \rightrightarrows M\otY \D_Y \rightarrow
M\ot_{\D_Y} \D_Y = M
\end{equation}
defining the $\D$-tensor product for any quasicoherent right $\cD_Y$-module
$M$.
The map $M\otY\D_Y\rightarrow
M$ is split surjective as a map of right $\OY$-modules (thanks to the
unit map $\OY\rightarrow\D_Y$). Taking $f^!$ of
\eqref{coequalizer defining tensor product},
then, gives a
surjective map
\begin{equation}\label{first arrow}
 f^!(M\otY\D_Y) = M\otY\DYX \rightarrow f^!(M);
\end{equation}
here the equality on the left follows from \eqref{shriek eq}.
The composites
of \eqref{first arrow} with the two arrows
\begin{equation}\label{double arrow for DYX}
f^!(M\otY\D_Y\otY\D_Y) = M\otY\D_Y\otY\DYX \rightrightarrows M\otY\DYX
\end{equation}
agree by construction.
 This implies that the coequalizer of \eqref{double arrow for DYX}, which is
$M\ot_{\D_Y}\DYX$, must map surjectively to $f^!(M)$.  In
particular, by Lemma \ref{lemma from mel}, $M\ot_{\D_Y}\DYX = 0$
implies $M = 0$.

The conclusion of the last paragraph applied to $M=\coker(\Phi)$
 thus implies that $\Phi$ is surjective.
It then follows that $\Phi$ is injective as well; for, letting
$K=\ker(\Phi)$, applying $f^*$  to the exact sequence
\bd
0\rightarrow K\rightarrow \DYX\otimes_{\D_X} \DXY \xrightarrow{\Phi} \D_Y
\rightarrow 0
\ed
gives a sequence
\bd f^*K \rightarrow
\DXY \rightarrow \DXY\rightarrow 0, \ed
in which the map from $f^*K$
must be injective because $\D_Y$ is flat over $\OY$ (Corollary
\ref{flatness cor}).  But then $f^*K =
0$, which implies that the surjective map \bd \DYX\ot_{\D_X}f^*K =
\DYX\ot_{\D_X}\DXY\ot_{\D_Y}K \rightarrow K \ed must have zero image,
i.e. $K=0$.  This completes the proof of the theorem.
\end{proof}

Combining this result with Proposition \ref{crystals descend} we
obtain the following:
\begin{corollary}
Suppose $f:X\rightarrow Y$ is a cuspidal quotient
morphism from  a smooth $k$-variety $X$.  Then:
\begin{enumerate}
\item The bimodules $\DXY$ and $\DYX$ induce Morita equivalences of
the categories of (left and right) $\D_X$-modules and $\D_Y$-modules.
\item
The category of right $\D_Y$-modules is
 equivalent to the category of $!$-crystals on $Y$ and thus
is equivalent to the category of
right $\D_Z$-modules supported on $Y$ for any
closed embedding $Y\hookrightarrow Z$ of $Y$ in a smooth variety $Z$.
\end{enumerate}
\end{corollary}

\begin{remark}\label{the D bimodules}
The bimodules $\DXY$ and $\DYX$ (Definition \ref{def of bimodules})
are adaptations to the singular context of the standard bimodules used
for pushforward and pullback of $\D$-modules. The standard definition
of the bimodule $\DXY$ (used for pullback of left $\D$-modules) is
$\DXY=(f\times 1)^*\D_Y$, which agrees with our definition thanks to
Lemma \ref{pullbacks of D}.  Thus for any left $\D_Y$-module, we have
a canonical isomorphism of $\D_X$-modules
$f^*M=\DXY\underset{\D_Y}{\ot}M$. Note in particular that under our
equivalence, the left $\D_Y$-module $\OY$ corresponds to the left
$\D_X$-module $\OX$.

The definition given for $\DYX$ is
likewise the replacement, for the not necessarily Gorenstein variety
$Y$, of the standard definition of the bimodule $\DYX$ using the
canonical line bundle.  Recall that for a morphism of {\em smooth}
varieties $f:X\to Y$, one defines an $(f\inv\D_Y,\D_X)$-bimodule
$\DYX$ as $\omega_X\otX f^*(\D_Y\otY\omega_Y\inv)$, where we first
turn the right $\D_Y$-structure of $\D_Y$ into a left structure, pull
back to a left $\D_X$-module, and then reconvert to a right
$\D_X$-module, using the canonical sheaves of $Y$ and $X$.  Thus, for
a smooth morphism $f$ we may write $\DYX=f^*\D_Y\otX \omega_{X/Y}$,
the tensor product with the relative canonical bundle (here pullback
is along the right structure---otherwise one utilizes the
transposition isomorphism \cite{Be2}). However, for a smooth affine
morphism we also have the identification of $f^!M$ with the
$\OX$-module corresponding to the $f_*\OX$-module
$\uHom_{\OY}(f_*\OX,M)=\uHom_{\OY}(f_*\OX,\OY)\ot M$, namely
$f^!M=f^*M\ot\omega_{X/Y}$. Thus, in this case $\DYX=(1\times
f)^!\D_Y$ and our definitions agree (again by Lemma \ref{pullbacks of
D}).

Note also that for any quasicoherent right $\D_Y$-module $M$ we have
$f^!M = M\otimes_{\D_Y}\DYX$.  Indeed, the proof of
Theorem \ref{the Morita theorem} shows that the map
$M\otimes_{\D_Y}\DYX \rightarrow f^!M$ is surjective.  It then follows by
a diagram chase that if $M\rightarrow N$ is a surjective map of
quasicoherent right $\D_Y$-modules then $f^!M\rightarrow f^!N$ is
surjective, so $f^!$ is an exact functor on quasicoherent right
$\D_Y$-modules.  Since $f^!\D_Y^I = \D_Y^I\otimes_{\D_Y}\DYX$, taking
a presentation of $M$ we find that $f^!M = M\otimes_{\D_Y}\DYX$ as well.

The definition using jets or equivalently $!$-pullback has
obvious functorial advantages over the definition using the canonical
sheaf and, as we have seen, has the correct role in the cuspidal
setting. Thus the Morita equivalences of Theorem \ref{the Morita
theorem} are given by the suitable adaptations of the pullback and
pushforward functors for $\D$-modules.
\end{remark}

\subsection{Dripping Varieties}\label{dripping varieties section}
\begin{defn}
The {\em deRham space} $Y_{dR}$ of a variety $Y$ is the quotient, in
the category of spaces (that is, {\it fppf} sheaves of sets), of $Y$ by the
formal groupoid $\wh{Y\times Y}$.
\end{defn}
As a functor, $Y_{dR}$ assigns to a scheme $S$ the set $Y_{dR}(S)=Y(S^{red})$, 
the set of
$Y$-points of the reduced scheme of $S$.  As a
$k$-space, $Y_{dR}$ is far from being algebraic.
\begin{defn}
A ($*$-)coherent sheaf on $Y_{dR}$ is an equivariant sheaf for the
formal groupoid of the diagonal under $*$-pullback, in other words a comodule 
for
$\cJ_Y$. We can also define a $!$-coherent sheaf on $Y_{dR}$ as an equivariant 
sheaf under
$!$-pullback, namely a $\cJ_Y$-cocomodule.
\end{defn}

For $f:X\to Y$ a cuspidal quotient, the deRham spaces $X_{dR}$ and
$Y_{dR}$ are very similar: the groupoid $X\times_Y X$ defining the
quotient map $X\to Y$ is a subgroupoid of the deRham groupoid of $X$,
so we may expect the map $X\to X_{dR}$ to factor through $X\to Y$ and
to identify the deRham spaces of $X$ and $Y$. The directed system
of all cusp quotients of $X$ can be considered a finitary
approximation of the deRham quotient. Thus we successively quotient
out by increasing finite infinitesimal equivalence relations in order
to approach the quotient of $X$ by the full infinitesimal nearness
relation $\wh{X\times X}$: the smooth variety $X$ ``drips'' down
towards $X_{dR}$ through the cuspidal varieties $Y$.

\begin{corollary} Suppose that
$X\rightarrow Y$ is a good cuspidal quotient morphism
between good CM varieties $X$ and $Y$.  Then the
deRham spaces $X_{dR}$ and $Y_{dR}$ have equivalent categories of
$*$-coherent sheaves and $!$-coherent sheaves (respectively).
\end{corollary}

Thus the dripping varieties picture is accurate on the level of
coherent sheaves: the ($*$ or $!$) pullback of $\OY$-modules from
deepening cusps to $X$ gives rise to sheaves with an increasingly
large piece of a (left or right) $\D_X$-module structure, while the
pullback of $\D$-modules from any cusp to $X$ gives rise to all
$\D_X$-modules.

\section{Cusp Induction}\label{cusp induction section}
Throughout this section, $X\xrightarrow{f} Y$ will always denote a good
cuspidal quotient morphism between good Cohen-Macaulay varieties over
the field $k$, and all $\D$-modules and $\theo$-modules
 are assumed to be quasicoherent.
\begin{notation}
Let $\pi_{X\times X},\pi_{Y\times X}$, and $\pi_{Y\times Y}$ be the
projections onto the first factor from the formal completions of
$X\times X$ along the diagonal, $Y\times X$ along the graph of $f$, and
$Y\times Y$ along the diagonal, respectively.
\end{notation}
We recall the notion of induced $\D$-module from \cite{S} (see also
\cite{chiral}). There is an exact faithful functor
$\Ind_X:\OX\mbox{--}\on{mod}\to\on{mod}\mbox{--}\D_X$ that sends an
$\OX$-module $M$ to the induced $\D_X$-module
$$\on{Ind}_X(M)=M\otX\D_X=\uHom_{\on{pro-mod}(\OX)}(\cJ_X,M)=\pi_{X\times
X}^!M.$$ The functor naturally lands in $(\OX,\D_X)$-bimodules, and we
then forget the commuting $\OX$-structure.

One may define a similar induction functor from $\OY$-modules to
$\cJ_Y$-cocomodules for any scheme $Y$.  We wish to construct a large
category of $\D_X$-modules on $X$, sharing some of the good properties
of induced $\D_X$-modules, by collecting induced modules from all
cuspidal quotients of $X$.

\subsection{Exactness of Cusp Induction}
As before, let $\JYX$ denote the pro-coherent $(\OY,\OX)$-bimodule (and 
commutative
algebra) of functions on the formal completion $\wh{Y\times X}$ of
$Y\times X$ along the graph of $f:X\to Y$.

\begin{prop}\label{projectivity of jets}\mbox{}
 \begin{enumerate}
\item The functor
$$M_Y\mapsto\pi_{Y\times X}^!M_Y=\uHom_{\OY}(\JYX,M_Y)$$ on
$\OY\mbox{--}\on{mod}$ is exact and equivalent to $M_Y\mapsto M_Y\otY\DYX$.
\item The functor $\pi_{Y\times X}^!$ may be refined to an exact
functor from $\OY\mbox{--}\on{mod}$ to $\rmod\mbox{--}\D_X$, the functor of 
{\em cusp
induction} $M_Y\mapsto \Ind_{Y\leftarrow X} M_Y$.
\end{enumerate}
\end{prop}

\begin{proof}
Part (1) is proven in Section \ref{section on jets and D-modules}.  For part 
(2), the functor
$M_Y\mapsto M_Y\otY \DYX$ clearly lands in
the category of right $\D_X$-modules, and since the
underlying functor to $\OX$-modules is exact, it follows that the
refined functor is exact as well.
\end{proof}

\subsection{Cusp-Induced Modules}\label{cusp-induced modules section}
Let $M_Y$ denote an $\Oo_Y$-module, and $M_X=f^*M_Y$ the corresponding
$\OX$-module.  Let $\Ind_X M_X=M_X\otX\D_X$ and
$\Ind_Y M_Y=M_Y\otY \D_Y$ be the corresponding induced right $\D_X$-module and
right $\D_Y$-module. The cusp-induced $\D_X$-module
$\Ind_{Y\leftarrow X} M_Y =M_Y\otY\DYX$ is an intermediate object.
\begin{lemma}
We have canonical isomorphisms
\begin{enumerate}
\item $(f\times 1)^*\Ind_{Y\leftarrow X}M_Y=\Ind_X M_X$ with respect to the 
left $\OY$-structure on
$\Ind_{Y\leftarrow X}M_Y$, and
\item $(1\times f)^! \Ind_Y M_Y =\Ind_{Y\leftarrow X}M_Y$ with
respect to the right $\OY$-structure on $\Ind_Y M_Y$.
\end{enumerate}
\end{lemma}
\noindent
This follows immediately from Lemma \ref{pullbacks of D}.

Recall from Section \ref{cusp quotients and jets} the definition
of $\D_X(N,M)$.
\begin{notation}
If $M_Y$ is an $\OY$-module equipped with an $\OY$-embedding
$M_Y\subseteq M$ in an $\OX$-module $M$ and $N$ is an $\OX$-module, we
let \bd \D(N,M_Y)\subseteq \D_X(N, M) \ed denote the subsheaf that
consists of those operators the image of which lies in $M$.
\end{notation}
\begin{prop}  Suppose that $X$ is a smooth $k$-variety.
With the above notation, we have
\bd
\D_Y(\OX, M_Y) = \D(\OX,M_Y).
\ed
\end{prop}
\begin{proof}
Because $Y$ is a cuspidal quotient of a smooth $k$-variety, it is very good
Cohen-Macaulay; hence the functor $M_Y\mapsto \uHomY(\JYX, M_Y)$ is exact.
In particular, the embedding $M_Y\rightarrow M$ induces an embedding
\bd
\D_Y(\OX, M_Y)\hookrightarrow \uHomY(\JYX,M) = \uHomX(\cJ_X, M) = \D_X(\OX, M).
\ed

On the other hand, $\D_X(\OX, M)$ consists of those $k$-linear maps
$\theta:\OX\rightarrow M$ for which $I_{\Delta_X}^n\cdot\theta = 0$ for
some $n\geq 0$.  So if $\theta\in\D_X(\OX,M)$ takes $\OX$ into $M_Y$, we
have in particular $I_{\Delta_Y}^n\cdot\theta = 0$ and so $\theta$ lies
in $\D_Y(\OX,M_Y)$ as well, completing the proof.
\end{proof}

\begin{corollary}[Agreement
with the Cannings-Holland construction, \cite{CH ideals}]\label{CH induction}

If $M_Y$ is a rank 1 torsion-free $\OY$-module on a birational
cuspidal quotient $X\rightarrow Y$ of a nonsingular curve $X$ and
$M_Y\subset K_X$ is an embedding of $M_Y$ as an $\theo_Y$-module, then
the cusp-induced $\D_X$-module $\widetilde{M}$ equals the subsheaf
$\D_X(\OX,M_Y)\subset \D_X(\theo_X,K_X)$ of operators with image
in $M_Y$.
\end{corollary}

\subsection{Riemann-Hilbert Correspondence for Cusp-Induced Modules}\label
{cusp RH section}
Recall that an {\em $f$-differential}
quasicoherent $(\OY,\OX)$-bimodule $M$ is a quasicoherent
sheaf on $Y\times X$ that
is torsion with respect to the ideal $I_{Y\leftarrow X}$ of the
graph of $f$ in $Y\times X$ (i.e., it is the union of its sections
supported on finite-order neighborhoods of the graph of $f$, and hence
underlies a $\JYX$-module).  We consider $M$ as a sheaf of
$(\OY,\OX)$-bimodules on $Y$.  Conversely \cite[Section 1.1.3]{BB},
if $M$ is a sheaf of $(\OY,\OX)$-bimodules on $Y$ that is quasicoherent
as an $\OY$-module and is torsion with respect to the ideal
$I_{Y\leftarrow X}\subset \OY\otimes_k\OX$ of the graph of $f$,
then $M$ is $f$-differential.

An $(\OY,\D_X)$-bimodule is said to be
$f$-differential if the underlying $(\OY,\OX)$-bimodule is.  The
category of $f$-differential $(\OY,\D_X)$-bimodules is denoted
$\on{mod}\mbox{--}(\OY,\D_X)_f$.

The deRham cohomology of right $\D_X$-modules $h_X:M\to
M\ot_{\D_X}\OX$ defines a functor from differential bimodules to
$\OY$-modules,
$$h_X:\on{mod}\mbox{--}(\OY,\D_X)_f\longrightarrow\OY\mbox{--}\on{mod}.$$
Note that we do not need to sheafify $h_X$ on the category of
differential bimodules, since $M\ot_{\D_X}\OX$ is automatically a
quasicoherent $\OY$-module.

\begin{thm} Cusp induction defines a fully faithful functor
$$\Ind_{Y\leftarrow X}:\OY\mbox{--}\on{mod}\longrightarrow
\on{mod}\mbox{--}(\OY,\D_X)_f,$$ which is right adjoint and right
inverse to the deRham functor $h_X$.
\end{thm}
\begin{proof}
A cusp-induced $\D_X$-module $\Ind_{Y\leftarrow X}(M_Y) =M_Y\otY\DYX$
carries an $f$-local $\OY$-structure commuting with its right
$\D_X$-structure, since it is the tensor product of a coherent sheaf
with the $f$-differential module $\DYX$; thus we may refine cusp
induction to a functor to $(\OY,\D_X)_f$-bimodules. Furthermore, the
deRham functor takes $(\OY,\D_X)$-bimodules to $\OY$-modules.

By Theorem \ref{the Morita theorem}, $\DYX$ is flat as a right $\D_X$-module, 
and
$\DYX\stackrel{\bbL}{\underset{\D_X}{\ot}}\OX=\OY$. It follows that for any 
$\OY$-module $M$,
$$\Ind_{Y\leftarrow X}M_Y\stackrel{\bbL}{\underset{\D_X}{\ot}}\OX=
M_Y\otY\DYX\stackrel{\bbL}{\underset{\D_X}{\ot}}\OX=M_Y,$$ so that
cusp-induced $\D$-modules are deRham exact, and the deRham functor
applied to an induced $(\OY,\D_X)$-bimodule $\Ind_{Y\leftarrow X}M_Y$
recovers $M_Y$ as $\OY$-module.

If $\Mtil$ is an $(\OY,\D_X)$-bimodule that is isomorphic to an
induced bimodule $\wt{N}= N_Y\otY\DYX$, then we have $h(\Mtil)\cong
h(\wt{N})= N_Y$, so $\Mtil$ is isomorphic to the induction of its
deRham module.  It follows that the deRham and induction functors
define quasi-inverse equivalences of categories between the essential
image of induction in $(\OY,\D_X)$-bimodules and $\OY$-modules,
provided we can show that the induction functor is full, which is a
special case of the adjunction.

Suppose $M$ is an $f$-differential bimodule and $N$ is a quasicoherent
$\OY$-module. We would like to show that $\OY$-morphisms $h_X(M)\to N$
are in bijection with $(\OY,\D_X)$-morphisms $M\to
\uHom_{\OY}(\JYX,N)$. Consider the diagram \bd \xymatrix@C+3ex{
\uHom_{\OX}(\cJ_X, M)\ar[d]^{(1)}
\ar@<.5ex>[r]^(.61){u_r^*}\ar@<-.5ex>[r]_(.61){a} & M\ar[d]^{(2)}
\ar[r]^(.6){p_M} & h_X(M)\ar@{-->}[d]^{(3)} \\
\uHom_{\OX}(\cJ_X,\uHom_{\OY}(\JYX,N))
\ar@<.5ex>[r]^(.61){u_r^*}\ar@<-.5ex>[r]_(.61){\Delta^*} &
\uHom_{\OY}(\JYX, N) \ar[r]^(.6){u_l^*} & N} \ed where $a$ and
$\Delta^*$ are the action maps, and the rows are colimit
diagrams. If we are given an $\OY$-module map $(3)$, the
composite $(3)\circ p_M$ induces, by adjunction, a $\cJ$-module map
$(2)$ such that $u_l^*\circ (2) = (3)\circ p_M$: the map is defined by
$(2)(m)(j_1\otimes j_2) = (3)\circ p_M(j_1\otimes j_2\cdot m)$ for
$j_1\otimes j_2$ a section of $\cJ$ and $m$ a section of $M$.  If
$(1)$ is the map induced by $(2)$ under $!$-pullback, then $(2)\circ
u_r^* = u_r^*\circ (1)$, and consequently $(2)\circ a = \Delta^*\circ
(1)$: we have \bd (3)\circ p_M\circ a = (3)\circ p_M \circ u^*_r =
u_l^*\circ (2)\circ u_r^* = u_l^*\circ u_r^*\circ (1)= u_l^*\circ
\Delta^*\circ(1), \ed and so for a section $s$ of $\uHom_Y(\cJ,M)$ we
find that
\begin{multline*}
(2)(a(s))(j_1\otimes j_2) = (3)\circ p_M (j_1\otimes j_2\cdot a(s)) =
(3)\circ p_M \circ a(j_1\otimes 1\otimes j_2\cdot s)\\ =
u_l^*\Delta^*(1)(j_1\otimes 1\otimes j_2\cdot s) = u_l^*(j_1\otimes
j_2\cdot \Delta^*(1)(s)) = \Delta^*\circ(1)(s)(j_1\otimes j_2).
\end{multline*}
Thus $(2)$ is an $(\OY,\cD_X)$-bimodule map.  Conversely, starting
from an $(\OY,\cD_X)$-bimodule map $(2)$, one has $(2)\circ a =
\Delta^*\circ (1)$ and $(2)\circ u_r^* = u_r^*\circ (1)$, so the
composite $u_l^*\circ (2)$ factors through $h(M)$ as $(3)\circ p_M$
for some $\OY$-module map $(3)$.  Under adjunction, moreover,
$(3)\circ p_M$ corresponds to $(2)$, as desired.  This establishes the
adjunction of $h_X$ and $\Ind_{Y\leftarrow X}$, completing the proof.
\end{proof}

\begin{lemma}\label{inductions agree}
Let $g:Y\to Y'$ be a morphism of cusp quotients of $X$.
Then the functors $\Ind_Y$ and $\Ind_{Y'}\circ g_*$ from
$\OY$-modules to $\D_X$-modules are naturally isomorphic.
\end{lemma}

\begin{proof}
The lemma follows from the identification $(g\times 1)^*\cJ_{Y'\leftarrow X}
=\JYX$
(under the morphism $g\times 1:Y\times X\to Y'\times X$)
and the ensuing identification
$$\Ind_{Y'\leftarrow X}g_*M_Y=\uHom_{\Oo_{Y'}}(\cJ_{Y'\leftarrow X},g_*M_Y)
=\uHom_{\Oo_{Y}}(g^*\cJ_{Y'\leftarrow X},M_Y)=\Ind_{Y\leftarrow X} M_Y,$$ 
which respects
cocomodule structures.
\end{proof}

\begin{defn}
For $\OY$-modules $M_Y$ and $N_Y$, the {\em differential morphisms}
are the $\D_X$-module homomorphisms between
cusp-induced modules,
$$\on{Diff}(M_Y, N_Y)=\on{Hom}_{\rmod-\D_X}(\on{Ind}_Y M_Y,\on{Ind}_Y N_Y).$$
\end{defn}

\begin{defn}
The category $\on{cusp-ind}(\cD_X)$ of {\em cusp-induced
$\cD$-modules} on $X$ is the full subcategory of $\on{mod}\mbox{-}\D_X$
consisting of $\D_X$-modules isomorphic to $\Ind_{Y\leftarrow X} M_Y$ for some 
cusp
quotient $X\xrightarrow{f} Y$ and some $\OY$-module $M_Y$.
\end{defn}

\begin{corollary}[{\bf Cuspidal Riemann-Hilbert Correspondence}]
\label{cusp RH}
The functors $\Ind_{Y\leftarrow X}$ define an equivalence of categories
\bd \underset{\underset{X\rightarrow
Y}{\longrightarrow}}{\lim} \left(\operatorname{qcoh}(\theo_Y),
\operatorname{Diff}\right) \longrightarrow
\operatorname{cusp-ind}(\cD_X).  \ed A quasi-inverse functor is given
by the deRham functor.
\end{corollary}

\begin{proof}
The functor $\Ind_{Y\leftarrow X}$ from $\OY\mbox{--}\on{mod}$ to $\on{mod}
\mbox{--}\D_X$ is faithful,
since the $(\OY,\D_X)$-morphisms between cusp-induced modules inject
into the $\D_X$-morphisms. It follows that by allowing differential
morphisms we make this faithful functor full. By Lemma \ref{inductions agree}, 
it is compatible with
deepening the cusps and so descends to the inductive limit. Since cusp-induced
$\D_X$-modules are by definition the essential image of this functor,
it follows that we have an equivalence of categories.
\end{proof}

The corollary is an extension of the equivalence between induced
$\D_X$-modules and $\OX$-modules with differential morphisms
\cite{S}.

\subsection{Cusp-Induced $\cD$-Modules on a Curve}\label{cusp-induced on curve 
section}

Let $X$ denote a smooth curve over a field $k$
of characteristic zero; the category of
cusp-induced $\cD$-modules on $X$ is particularly easy to describe in
this case. The cusp-induction functor itself also becomes very
concrete: as we have seen (Corollary \ref{CH induction}), in the case
of curves, the cusp-induction of an $\OY$-module $M\subset K_Y=K_X$
agrees with the Cannings-Holland construction, the sheaf $\D_X(\OX,M)$
of differential operators with values in $M$.

\begin{defn}
For a right $\cD_X$-module $M$, we let $T_{\theo}(M)$ denote the
$\OX$-torsion submodule of $M$, and $T_{\cD}(M)$ denote the
$\cD_X$-torsion submodule of $M$.  We say $M$ is {\em generically
torsion-free} if there is a nonempty open set $U$ of $X$ such that
$M|_U$ is a torsion-free $\cD_U$-module, in other words if
$T_{\cD}(M)$ is supported on a proper closed subvariety.
\end{defn}

\begin{lemma}
$T_\theo(M)$ and $T_\cD(M)$ are $\cD_X$-submodules of $M$, which agree if $M$
is generically torsion-free.
\end{lemma}

\begin{prop}\label{cusp-induced on a curve}
A finitely-generated $\cD_X$-module $M$ is cusp-induced if and only if
$M$ is generically torsion-free.
\end{prop}
\begin{proof}
Suppose $M = \Mbar\otY\DYX$ for a cuspidal quotient $X\rightarrow Y$.
$\Mbar$ is an extension \bd 0\rightarrow T_\OY(\Mbar) \rightarrow
\Mbar \rightarrow \Mbar/T_\OY(\Mbar) \rightarrow 0 \ed where
$T_\OY(\Mbar)$ is $\OY$-torsion and $\Mbar/T_\OY(\Mbar)$ is
$\OY$-torsion-free.  Since $\DYX$ is flat over $\OY$, we obtain an
exact sequence \bd 0\rightarrow T_\OY(\Mbar)\otY\DYX \rightarrow M
\rightarrow (\Mbar/T_\OY(\Mbar))\otY \DYX \rightarrow 0, \ed where
$T_\OY(\Mbar)\otY\DYX$ is an $\OX$-torsion submodule of $M$.
Moreover, because $\Mbar/T_\OY(\Mbar)$ is $\OY$-torsion-free, it
embeds locally in $\Oo_Y^n$ for some $n$, hence
$(\Mbar/T_\OY(\Mbar))\otY\DYX$ embeds locally in $\cD_X^n$, so it is
$\cD_X$-torsion-free.

Conversely, suppose $M$ is a $\cD_X$-module which is generically torsion-free.
The $\OX$-torsion submodule $T_\theo(M)$ is supported on a finite
collection of closed points, hence by Kashiwara's Theorem it is
isomorphic to an induced $\cD_X$-module.  Moreover, $M/T_\theo(M)$ is
$\cD_X$-torsion-free, hence is locally projective on $X$; therefore, the group
$\Ext_{\cD_X}^1(M/T_\theo(M), T_\theo(M))$ vanishes and so $M \cong
T_\theo(M)\oplus M/T_\theo(M)$.  Now $M/T_{\theo}(M)$ embeds in a locally free
(and therefore induced) $\cD_X$-module with cosupport a finite set.
It thus suffice to prove the following easy consequence
of Kashiwara's theorem.
\begin{lemma}\label{finite cosupport means cusp-induced}
Suppose
\bd
{\mathbf L}: \;\;\;
0\rightarrow \ker(\beta)\rightarrow M\xrightarrow{\beta} Q\rightarrow 0
\ed
is a short exact sequence of finitely generated $\D_X$-modules such that
$M$ is cusp-induced and $Q$ is supported on a finite subset of $X$.
 Then there is a cuspidal quotient $X\rightarrow Y$ and an exact sequence of
coherent $\OY$-modules
\bd
{\mathbf L'}:\;\;\;
0\rightarrow K\rightarrow M' \rightarrow Q'\rightarrow 0
\ed
such that $\Ind_{Y\leftarrow X}({\mathbf L'}) \cong {\mathbf L}$ as sequences
of right $\D_X$-modules.
\end{lemma}
\noindent
This completes the proof of Proposition \ref{cusp-induced on a curve}.
\end{proof}

The sheaf of algebras $\D_X$ locally has homological dimension
one (see \cite[Section 1.4(e)]{SS}), from which it follows that any torsion-
free $\D$-module is
locally projective. Moreover, $\D_X$ possess a skew field of fractions
(see \cite[Section 2.3]{SS}), from
which it follows that any finitely generated
locally projective $\D$-module has a well-defined
rank.

\begin{defn}\cite{chiral} A $\D$-vector bundle (or $\D$-bundle for short)
$M$ on $X$ is a locally projective (equivalently torsion-free) right
$\D_X$-module of finite rank.
\end{defn}

\begin{corollary} A $\D$-module $M$ on $X$ is a $\D$-bundle if and only if
it is isomorphic to a $\D$-module cusp-induced from a torsion-free
$\OY$-module for some cuspidal quotient $Y$ of $X$.
\end{corollary}

An important class of $\D$-bundles (of rank one) is provided by the right 
ideals in
$\D_X$ (in fact any rank one $\D$-bundle may locally
 be embedded as a right ideal in $\D_X$).

\begin{example}
Consider the right ideal $M\subset\D_{\aline}=\C\langle
z,\del\rangle/\{\del z-z\del-1\}$ generated by $z^2$ and $1-z\del$.
Then $M=\D(\theo_{\aline},\OY)$ where $\OY=k[z^2,z^3]$ is the
coordinate ring of the
cuspidal cubic curve $y^2=x^3$, i.e. $M$ is cusp-induced from the
structure sheaf of $Y$.
\end{example}

As the example illustrates, $\D$-bundles are not locally trivial---that is,
they are not locally (on the base curve)
 isomorphic to a direct sum of copies of $\D$. It follows from the
cusp-induced description, however, that $\D$-bundles are {\em
generically} trivial.  It is convenient to parametrize $\D$-bundles by
picking such a generic trivialization. This leads to the description
by Cannings and Holland \cite{CH ideals, CH
cusps}
of ideals in $\D$
by means of the {\em adelic Grassmannian}
\cite{Wilson}.

\begin{defn} The adelic Grassmannian $\Gr^{ad}(X)$ is the set of isomorphism 
classes
of $\D$-bundles $M$ of rank one equipped with a generic
trivialization, $M\ot K_X\cong \D_X\otimes K_X$.
\end{defn}
\begin{remark}
This set-theoretic definition may be refined to a moduli problem,
giving the adelic Grassmannian an algebraic structure, which we study
in \cite{W}. More precisely, for a finite set $I$, there is an
ind--scheme of ind--finite type $\Gr^{ad}_I(X)$ over $X^I$,
parametrizing $\D$--bundles with trivializations outside of
$I$--tuples of points of $X$.  As we allow these prescribed positions
of singularities to collide, we obtain a ``nonlinear vertex algebra''
structure on $\Gr^{ad}(X)$: the spaces $\{\Gr^{ad}_I(X)\}$, considered
over the inductive system of spaces $X^I$ with respect to partial
diagonal maps (i.e. surjections of finite sets), form a factorization
ind--scheme, as defined in \cite{chiral}. We also identify the vertex
(or chiral) algebras obtained by linearization of this factorization
space with the ${\mathcal W}_{1+\infty}$-vertex algebra and its
variants.
\end{remark}

\begin{corollary} $\Gr^{ad}(X)$ is isomorphic to the direct limit over cusp 
quotients
$X\to Y$ of the set of isomorphism classes of rank $1$ torsion-free
$\OY$-modules equipped with a generic trivialization.
\end{corollary}

\begin{proof}
Observe that if $X\rightarrow Y$ is a cuspidal quotient morphism (over
 $k$ of characteristic zero) then it is birational.  It follows that
 $K_Y\otimes_{\OY}\DYX = K_X\otimes_{\OX}\D_X$ canonically, so by
 Corollary \ref{cusp RH} we get an embedding of the set of isomorphism
 classes of rank $1$ torsion-free $\OY$-modules equipped with generic
 trivialization in $\Gr^{ad}(X)$. (Note that since we have rigidified
 our modules using the generic trivialization, the distinction between
 differential morphisms and $\OY$--morphisms goes away -- this
 injectivity is certainly false without the rigidification.) It
 suffices, then, to prove that this embedding is surjective.  So
 suppose $M$ is a rank one $\D$-bundle equipped with an isomorphism
 $\phi: M\otimes K_X \rightarrow \D\otimes K_X$.  Then there is an
 effective divisor $D$ on $X$ such that $\phi$ is of the form
 $\psi\otimes K_X$ for a $\D$-module homomorphism $\psi: M\rightarrow
 \theo(D)\ot_{\OX}\D_X$ that is an isomorphism on a nonempty open set
 of $X$.  Now Lemma
\ref{finite cosupport means cusp-induced} implies that there is a
cuspidal quotient $X\rightarrow Y$ and a
 homomorphism
$\psi': M'\rightarrow \theo(D)$ of $\OY$-modules such that
$\Ind_{Y\leftarrow X}(\psi') = \psi$.  The canonical map $\theo(D)\rightarrow 
K_X$ now gives the desired generic trivialization, completing the proof.
\end{proof}

As we let the cusps $Y$ get deeper, the rings $\OY$ evaporate so that
the adelic Grassmannian parametrizes simply certain linear algebra
data. This is most succinctly explained in
\cite{chiral}, Section 2.1: $\D$-submodules of any $\D$-module $M$ that
are cosupported at a point $x\in X$ are in canonical bijection (via
the deRham functor $h$) with subspaces of the stalk of the dRham
cohomology $h(M)_x$ at $x$ that are open in a natural topology. In the
above case, we have $\D$-submodules of $\D(K_X)$ cosupported at some
finite collection of points, which correspond to collections of open
subspaces of Laurent series ${\mathcal K}_{x_i}$ at $x_i$ with respect
to the usual topology.

\begin{remark}
The adelic Grassmannian (for $X=\aline$) was used by G. Wilson
to parametrize solutions of the KP hierarchy that are rational (and
decay at infinity) in the first time variable. More precisely,
these solutions correspond to Krichever data, which are rank one
torsion-free sheaves on a cuspidal quotient of the affine line. (More
generally $\Gr^{ad}(X)$ for a curve $X$ parametrizes Krichever data
for $X$ and all of its cuspidal quotients together.) Wilson then shows
that under the action of the KP flows on the adelic Grassmannian,
the underlying collections of points on $\aline$ move according to the
Calogero-Moser particle system. This extends the results of
Airault-McKean-Moser \cite{AMM}, Krichever \cite{Kr1} and Shiota
\cite{Sh} to allow collisions of the Calogero-Moser
particles. In \cite{solitons} we develop a $\D$-bundle point of view
on the KP hierarchy, and obtain geometric proofs of this result and
extend it to the rational, trigonometric and elliptic solutions of
(multicomponent) KP hierarchies.
\end{remark}

\end{document}